\documentclass[11pt]{amsart}
\usepackage{amssymb,amsmath,amsthm,mathtools,xcolor,booktabs}
\usepackage[margin=1.1in]{geometry}
\usepackage[colorlinks=true,linkcolor=blue!55!black,citecolor=blue!55!black,urlcolor=blue!55!black]{hyperref}
\newtheorem{theorem}{Theorem}[section]
\newtheorem{proposition}[theorem]{Proposition}
\newtheorem{lemma}[theorem]{Lemma}
\newtheorem{corollary}[theorem]{Corollary}
\newtheorem{conjecture}[theorem]{Conjecture}
\newtheorem{question}[theorem]{Question}
\theoremstyle{remark}
\newtheorem{remark}[theorem]{Remark}

\newtheorem{introthm}{Theorem}

\newtheorem{introprop}[introthm]{Proposition}
\newtheorem{introconj}[introthm]{Conjecture}

\newcommand{\Q}{\mathbb{Q}}
\newcommand{\Z}{\mathbb{Z}}
\newcommand{\F}{\mathbb{F}}
\newcommand{\Ep}{\mathcal{E}_p}
\newcommand{\B}{\mathcal{B}}
\newcommand{\CL}{\mathrm{CL}}
\newcommand{\R}{\mathcal{R}}

\title[Unobstructed deformation rings in the level aspect]{Unobstructed deformation rings, residual ramification, and modular degrees in the level aspect}
\author{Bartu Bingol}
\address{School of Mathematical and Statistical Sciences, Arizona State University, Tempe, AZ 85287, USA}
\email{bbingol@asu.edu}  
\date{\today}
\subjclass[2020]{11F80 (primary); 11G05, 11N36, 11F33 (secondary)}
\keywords{Galois deformation rings, unobstructed deformation problems, modular degree, congruence primes, arithmetic statistics, elliptic curves, level aspect}

\begin{document}

\begin{abstract}
\sloppy Fix a prime $p \geq 5$. For elliptic curves $E/\Q$ with good reduction at $p$, we study how often the deformation ring $\R_{E,p}$ of Ray and Weston is unobstructed. Ordering all curves by height, we prove that the number of unobstructed curves of height at most $X$ is bounded by a constant multiple of $X^{5/6}(\log X)^{-2/(p-1)}$. Using theorems of Mazur and Yasuda, we show that the local conditions defining the Ray-Weston family force $E[p]$ to be absolutely irreducible. Within this family, curves whose mod-$p$ representation is unramified at a bad prime are always obstructed and have $p$ dividing the modular degree. On the complement, data for all $2{,}164{,}260$ isogeny classes of conductor at most $500{,}000$, together with an Euler characteristic heuristic, suggest that $p$ divides the modular degree with the Cohen-Lenstra frequency of Watkins, with no detectable dependence on the number of bad primes.
\end{abstract}

\maketitle

\section{Introduction}

Let $p \geq 5$ be prime and let $E/\Q$ be an elliptic curve with good reduction at $p$ and mod-$p$ representation $\bar\rho = \bar\rho_{E,p} \colon G_\Q \to \mathrm{GL}_2(\F_p)$. When $\bar\rho$ is absolutely irreducible, Ray and Weston \cite{RW} attach to $(E,p)$ a deformation ring $\R_{E,p}$, parametrizing lifts of $\bar\rho$ unramified outside $p$ and the bad primes of $E$, with determinant the cyclotomic character and a flat (ordinary, in the split case) condition at $p$, and no condition at the bad primes. Since the determinant is fixed, $\R_{E,p}$ is unobstructed exactly when $\R_{E,p} \cong \Z_p$.

For a fixed curve and varying $p$ the picture is well understood: Weston \cite{W1} proved that the unrestricted deformation ring of a newform of squarefree level is unobstructed for a density-one set of primes, and Ray and Weston \cite[Cor.4.7]{RW} proved that for $E$ without complex multiplication and with squarefree conductor, $\R_{E,p} \cong \Z_p$ for all but finitely many $p$. The horizontal question, with $p$ fixed and $E$ varying, is \cite[Question 1.1(2)]{RW}. Ray and Weston studied it on the family $\Ep$ of curves with (1) good reduction at $p$, (2) squarefree conductor $N$, (3) $\bar\rho$ absolutely irreducible, and (4) $\ell \not\equiv \pm 1 \pmod p$ for every $\ell \mid N$ \cite[Def.5.1]{RW}. On $\Ep$ they showed that $p \nmid m_E$, where $m_E$ is the modular degree, implies $\R_{E,p} \cong \Z_p$, and computed the relevant proportions to conductor $4000$.

This paper has three results all depending on the local invariant of the behavior of $\bar\rho$ at a prime $\ell$ of multiplicative reduction, governed by whether $\ell \equiv \pm 1 \pmod p$ and whether $\bar\rho$ is unramified at $\ell$ (equivalently $p \mid v_\ell(\Delta_E)$). Throughout, $\kappa = \kappa_p = \frac{2}{p-1}$ and $R = R_p$ is the set of primes $\ell \equiv \pm 1 \pmod p$.

\subsection*{A density bound over all curves}

Order elliptic curves over $\Q$ by naive height $H$ (Section \ref{sec:upper}); there are $\asymp X^{5/6}$ of height at most $X$.

\begin{introthm}[Theorem \ref{thm:upper}]\label{ithm:upper}
Fix $p \geq 5$. Then
$$\#\big\{E/\Q :\ H(E) \leq X,\ p \text{ good},\ \bar\rho_{E,p} \text{ abs.\ irreducible},\ \R_{E,p} \cong \Z_p\big\} \ll_p X^{5/6}(\log X)^{-2/(p-1)}.$$
In particular, the proportion of elliptic curves for which $\R_{E,p}$ is unobstructed tends to $0$, at least as fast as $(\log X)^{-2/(p-1)}$.
\end{introthm}

The arithmetic input (Proposition \ref{prop:obstruct}) is that a single prime of multiplicative reduction in $R$ produces a nonzero local invariant $H^0(G_\ell, \mathrm{Ad}^0\bar\rho(1))$, which by the Euler characteristic formula of Ray and Weston forces a nonzero tangent space. No condition at primes of additive reduction is needed. The analytic input is an elementary upper-bound sieve: at a prime $\ell \geq 5$, multiplicative reduction cuts out exactly $\ell - 1$ of the $\ell^2$ residue classes of Weierstrass coefficients, and avoiding this set at every $\ell \in R$, a set of primes of relative density $\kappa$, costs a factor $(\log X)^{-\kappa}$. The same bound holds for the height-ordered count of $\Ep$ (Corollary \ref{cor:Ep}). We do not prove a matching lower bound; see Conjecture \ref{iconj:height}.

\subsection*{Irreducibility on \texorpdfstring{$\Ep$}{E\_p}}

Ray and Weston observed computationally that conditions (1), (2), (4) appear to imply (3), and noted that they were not aware of a reason \cite[Remark 5.3]{RW}.

\begin{introthm}[Theorem \ref{thm:irred}]\label{ithm:irred}
Let $p \geq 5$, and let $E/\Q$ have good reduction at $p$, squarefree conductor $N$, and $\ell \not\equiv \pm 1 \pmod p$ for every $\ell \mid N$. Then $E[p]$ is absolutely irreducible. For $p \geq 11$ the congruence condition on $N$ is not needed.
\end{introthm}

The proof is short: a reducible $E[p]$ forces some curve isogenous to $E$ to have a rational point of order $p$ (Proposition \ref{prop:shape}); Mazur's theorem \cite{M} rules this out for $p \geq 11$, and a theorem of Yasuda \cite{Y12} rules it out for $p \in \{5,7\}$ under the congruence condition on the bad primes.

\subsection*{Residual ramification and the modular degree}

Let $L$ denote the set of curves in the family for which $\bar\rho$ is unramified at some $\ell \mid N$. On $\Ep$ we have the inclusions
$$\{E \in \Ep : p \nmid m_E\} \subseteq \{E \in \Ep : \R_{E,p} \cong \Z_p\} \subseteq \Ep \setminus L$$
(Corollary \ref{cor:sandwich}): the first is \cite[Thm.4.6]{RW}, and the second is again Proposition \ref{prop:obstruct}. The set $L$ is also contained in $\{p \mid m_E\}$.

\begin{introprop}[Proposition \ref{prop:lowering}]\label{iprop:lowering}
Let $E$ satisfy conditions (1), (2), (4), and suppose $\bar\rho_{E,p}$ is unramified at some $\ell \mid N$. Then $p \mid m_E$.
\end{introprop}

This divisibility is assembled from known results: Mazur's principle lowers the level of $\bar\rho$ to a divisor of $N/\ell$, which makes $p$ a congruence prime for $f_E$ in the sense of \cite{RW}; the same implication appears in the proof of \cite[Thm.4.6]{RW}, where it is attributed to \cite[Lemma 12]{W2}. Its statistical consequence is the exact identity
$$P(p \mid m_E) = P(L) + \big(1 - P(L)\big) P(p \mid m_E \mid E \notin L)$$
for any finite set of curves in the family. In the data, the first term accounts for the growth of $P(p \mid m_E)$ with the number of bad primes, while the conditional frequency shows no statistically detectable variation across the observed strata (Table \ref{tab:omega}). This suggests:

\begin{introconj}[Conjecture \ref{conj:CL}]\label{iconj:CL}
Ordering the isogeny classes satisfying (1), (2), (4) by conductor, $P(p \mid m_E \mid E \notin L) \to \CL_p := 1 - \prod_{i \geq 1}(1 - p^{-i})$, and the same limit holds within each stratum of fixed $\omega(N)$.
\end{introconj}

$\CL_p$ is the Cohen-Lenstra-type constant proposed by Watkins \cite{Wat} and used in \cite{RW}. The evidence is as follows. At conductor $500{,}000$, a chi-square test does not reject homogeneity across $\omega(N)$ for any of $p = 5, 7, 11, 13$. The pooled conditional frequencies agree with $\CL_p$ to within sampling error for $p = 5, 7$, and for $p = 11, 13$ they lie below $\CL_p$ by an amount that shrinks steadily across disjoint conductor windows (Section \ref{sec:congruence}). At prime conductor the residual term is provably absent (Proposition \ref{prop:primevac}), which explains why the constant of \cite{RW} is well calibrated there. In Section \ref{sec:heuristic} we give a heuristic for Conjecture \ref{iconj:CL}: in the Euler characteristic formula for $\R_{E,p}$, the bad primes contribute only through the local invariants of Lemma \ref{lem:local}, and these vanish exactly off $L$, so off $L$ the relevant Selmer group is balanced at every bad prime and a Cohen-Lenstra model predicts the constant $\CL_p$ regardless of $\omega(N)$.

\begin{introconj}[Conjecture \ref{conj:height}]\label{iconj:height}
The bound of Theorem \ref{ithm:upper} has the correct order: the number of $E/\Q$ of height at most $X$ with $\R_{E,p} \cong \Z_p$ is $\gg_p X^{5/6}(\log X)^{-2/(p-1)}$.
\end{introconj}

We separate what is proved from what is observed. Theorems \ref{ithm:upper} and \ref{ithm:irred} and Proposition \ref{iprop:lowering} are proved. The ordering by height in Theorem \ref{ithm:upper} and the ordering by conductor in our data are different ensembles, and we do not prove any transfer between them. The data come from Cremona's database \cite{E} and reproduce the published conductor-$4000$ table of \cite{RW} (Section \ref{sec:repro}).

\subsection*{Organization}
Section \ref{sec:local} sets up the deformation problem and the local obstruction criterion. Section \ref{sec:irred} proves Theorem \ref{ithm:irred}, and Section \ref{sec:upper} proves Theorem \ref{ithm:upper}. Section \ref{sec:congruence} treats the residual-ramification locus and the conditional statistics. Section \ref{sec:Lstats} records what the data say, and do not yet say, about the size of $L$ in the conductor ordering. Section \ref{sec:repro} describes the data and computations, and Appendix \ref{app:counts} lists raw counts.

\section{The deformation problem and local obstructions}\label{sec:local}

We follow \cite[\S\S 2-3]{RW}. Let $E/\Q$ have good reduction at $p \geq 5$, let $S$ be the set of primes of bad reduction, and assume $\bar\rho = \bar\rho_{E,p}$ is absolutely irreducible. The ring $\R_{E,p}$ represents deformations of $\bar\rho$ to $G_{S \cup \{p\}}$ with determinant the cyclotomic character $\chi$, satisfying the liftable local condition $\mathcal{C}_p$ of \cite[Prop.2.9]{RW} at $p$ (with tangent space $\mathcal{N}_p$), and with no condition at the primes of $S$. Its tangent space is the Selmer group
$$H^1_{\langle p \rangle}(G_{S\cup\{p\}}, \mathrm{Ad}^0\bar\rho) = \ker\Big(H^1(G_{S\cup\{p\}}, \mathrm{Ad}^0\bar\rho) \to H^1(G_p, \mathrm{Ad}^0\bar\rho)/\mathcal{N}_p\Big),$$
and $\R_{E,p} \cong \Z_p$ forces this group to vanish, since the maximal ideal of $\Z_p$ is generated by $p$.

Write $\mathrm{Ad}^0\bar\rho^* = \mathrm{Ad}^0\bar\rho(1)$, the Tate dual of $\mathrm{Ad}^0\bar\rho$, and let $\mathcal{N}_p^\perp \subseteq H^1(G_p, \mathrm{Ad}^0\bar\rho^*)$ be the annihilator of $\mathcal{N}_p$ under local Tate duality. The Selmer structure defining $\R_{E,p}$ imposes the full local condition $H^1(G_\ell, \mathrm{Ad}^0\bar\rho)$ at each $\ell \in S$; its annihilator is $0$. The dual Selmer group is therefore
$$H^1_{\langle p \rangle}(G_{S\cup\{p\}}, \mathrm{Ad}^0\bar\rho^*) = \ker\Big(H^1(G_{S\cup\{p\}}, \mathrm{Ad}^0\bar\rho^*) \to \bigoplus_{\ell \in S} H^1(G_\ell, \mathrm{Ad}^0\bar\rho^*) \oplus H^1(G_p, \mathrm{Ad}^0\bar\rho^*)/\mathcal{N}_p^\perp\Big),$$
consisting of the classes that vanish locally at every bad prime and lie in $\mathcal{N}_p^\perp$ at $p$. (At the archimedean place $H^1(\mathbb{R}, \cdot) = 0$, as $p$ is odd.) With these definitions, the Greenberg-Wiles Euler characteristic formula gives \cite[Prop.3.2]{RW}
\begin{equation}\label{eq:wiles}
\dim H^1_{\langle p \rangle}(G_{S\cup\{p\}}, \mathrm{Ad}^0\bar\rho) - \dim H^1_{\langle p \rangle}(G_{S\cup\{p\}}, \mathrm{Ad}^0\bar\rho^*) = \sum_{\ell \in S} \dim H^0\big(G_\ell, \mathrm{Ad}^0\bar\rho(1)\big).
\end{equation}

\begin{lemma}[Local criterion]\label{lem:local}
Let $\ell \neq p$ be a prime of multiplicative reduction for $E$. Then
$$H^0\big(G_\ell, \mathrm{Ad}^0\bar\rho(1)\big) \neq 0 \iff \ell \equiv \pm 1 \pmod p \ \text{ or } \ \bar\rho \text{ is unramified at } \ell,$$
and $\bar\rho$ is unramified at $\ell$ if and only if $p \mid v_\ell(\Delta_E)$, where $\Delta_E$ is the minimal discriminant.
\end{lemma}

\begin{proof}
Over $\Q_\ell$, $E$ is a quadratic unramified twist of the Tate curve $E_q$, $q \in \Q_\ell^\times$ with $v_\ell(q) = -v_\ell(j_E) = v_\ell(\Delta_E) > 0$. Hence $\bar\rho|_{G_\ell} \cong \lambda \otimes \bar\rho_q$ with $\lambda$ unramified of order dividing $2$, and $\bar\rho_q$ is an extension of the trivial character by the mod-$p$ cyclotomic character $\bar\chi$, whose class in $H^1(G_\ell, \mu_p) = \Q_\ell^\times/(\Q_\ell^\times)^p$ is that of $q$. The field cut out by $\bar\rho_q$ is $\Q_\ell(\mu_p, q^{1/p})$, which is unramified over $\Q_\ell$ if and only if $p \mid v_\ell(q)$; this proves the last assertion.

In a basis adapted to the filtration, the line $\mathfrak{n}^+ \subset \mathrm{Ad}^0\bar\rho$ of strictly upper triangular matrices is $G_\ell$-stable with character $\bar\chi$, so $\mathfrak{n}^+(1)$ has character $\bar\chi^2$, which is trivial on $G_\ell$ when $\ell^2 \equiv 1 \pmod p$. This gives the implication from $\ell \equiv \pm 1$.

Suppose $\ell \not\equiv \pm 1 \pmod p$, so $\bar\chi|_{G_\ell}$ and $\bar\chi^2|_{G_\ell}$ are nontrivial. Since $p \nmid \ell - 1$, every unit of $\Z_\ell$ is a $p$-th power, so $\Q_\ell^\times/(\Q_\ell^\times)^p \cong \Z/p$ via $v_\ell$; thus $\bar\rho_q$ is split if and only if $p \mid v_\ell(q)$, i.e.\ if and only if $\bar\rho$ is unramified at $\ell$. A nonzero vector of $H^0(G_\ell, \mathrm{Ad}^0\bar\rho(1))$ is a nonzero $G_\ell$-map $\phi\colon \bar\rho \to \bar\rho \otimes \bar\chi$. The Jordan-H\"older constituents of the source are $\lambda\bar\chi, \lambda$ and those of the target are $\lambda\bar\chi^2, \lambda\bar\chi$; as $\bar\chi$ and $\bar\chi^2$ are nontrivial, the only constituent they share is $\lambda\bar\chi$. So $\phi$ has rank one, with image the line $\lambda\bar\chi^2$ or with image $\cong \lambda\bar\chi$ occurring as a quotient of $\bar\rho$. The first is impossible, since $\lambda\bar\chi^2$ is not a constituent of $\bar\rho$; the second requires $\lambda\bar\chi$ to be a quotient of $\bar\rho$, which happens exactly when $\bar\rho|_{G_\ell}$ is split. Conversely, in the split case $\mathrm{Ad}^0\bar\rho(1)$ is the sum of the characters $\bar\chi^2, \bar\chi, \mathbf{1}$, and has invariants. The case $\ell \not\equiv \pm 1$ of the lemma is also \cite[Prop.4.5]{RW}.
\end{proof}

\begin{proposition}\label{prop:obstruct}
Let $E/\Q$ have good reduction at $p \geq 5$ with $\bar\rho_{E,p}$ absolutely irreducible. If $E$ has multiplicative reduction at some prime $\ell$ with $\ell \equiv \pm 1 \pmod p$, or at some prime $\ell$ at which $\bar\rho_{E,p}$ is unramified, then $\R_{E,p} \not\cong \Z_p$.
\end{proposition}

\begin{proof}
By Lemma \ref{lem:local} the right side of \eqref{eq:wiles} is positive, so the tangent space $H^1_{\langle p \rangle}(G_{S\cup\{p\}}, \mathrm{Ad}^0\bar\rho)$ is nonzero.
\end{proof}

Write $r_E$ for the congruence number of $E$ \cite{ARS}. By \cite[Thm.2.1]{ARS}, $m_E \mid r_E$, and $\mathrm{ord}_p(r_E) = \mathrm{ord}_p(m_E)$ when $\mathrm{ord}_p(N) \leq 1$; in particular $p \mid r_E \iff p \mid m_E$ when $p \nmid N$.

\begin{proposition}[Ray-Weston]\label{prop:RW}
If $E \in \Ep$ and $p \nmid m_E$, then $\R_{E,p} \cong \Z_p$.
\end{proposition}

\begin{proof}
For $E \in \Ep$, the prime $p$ lies outside the exceptional set $\Sigma(E)$ of \cite[Thm.4.6]{RW} as soon as $p$ is not a congruence prime for $f_E$ in the sense of \cite[Def.4.3]{RW}. As noted in \cite[\S 5.1]{RW}, this holds exactly when $p \nmid r_E$, i.e.\ when $p \nmid m_E$.
\end{proof}

Let $L$ be the set of curves $E$ satisfying conditions (1), (2), (4) for which $\bar\rho_{E,p}$ is unramified at some $\ell \mid N$, equivalently (Lemma \ref{lem:local}) $p \mid v_\ell(\Delta_E)$ for some $\ell \mid N$. By Theorem \ref{thm:irred} below, such curves lie in $\Ep$. Since no curve in $\Ep$ admits a rational $p$-isogeny, all curves in an isogeny class in $\Ep$ have isomorphic mod-$p$ representations, so membership in $L$ is a property of the isogeny class.

\begin{corollary}\label{cor:sandwich}
Within $\Ep$,
$$\{p \nmid m_E\} \subseteq \{\R_{E,p} \cong \Z_p\} \subseteq \Ep \setminus L.$$
\end{corollary}

\begin{proof}
The first inclusion is Proposition \ref{prop:RW}, the second is Proposition \ref{prop:obstruct}.
\end{proof}

\begin{remark}[The converse on $\Ep \setminus L$]\label{rem:converse}
By Proposition \ref{prop:lowering} below, $L \subseteq \{p \mid m_E\}$, so the two outer sets in Corollary \ref{cor:sandwich} differ only by curves in $\Ep \setminus L$ with $p \mid m_E$. For such curves $p$ is a congruence prime, necessarily a strict one (a proper congruence would lower the level and place $E$ in $L$), and \cite[Prop.4.4]{RW} identifies this with nonvanishing of the Selmer group $H^1_\emptyset(\Q, W_p(E))$ of \cite[\S 4]{RW}. Equality in the \emph{first} inclusion of Corollary \ref{cor:sandwich}, i.e.\ $\{p \nmid m_E\} = \{\R_{E,p} \cong \Z_p\}$ on $\Ep$, would therefore follow from showing that $H^1_\emptyset(\Q, W_p(E)) \neq 0$ forces $H^1_{\langle p \rangle}(G_{S\cup\{p\}}, \mathrm{Ad}^0\bar\rho) \neq 0$, a comparison of local conditions at $p$ that we do not attempt here. The second inclusion is far from an equality: the data contain many classes outside $L$ with $p \mid m_E$ (Table \ref{tab:omega}). Accordingly, all statistics below are statistics of $p \mid m_E$; for any finite set of curves in $\Ep$ they give
$$1 - P(p \mid m_E) \leq P(\R_{E,p} \cong \Z_p) \leq 1 - P(L).$$
\end{remark}

\section{Irreducibility of \texorpdfstring{$E[p]$}{E[p]}}\label{sec:irred}

\begin{proposition}\label{prop:shape}
Let $E/\Q$ be semistable with good reduction at $p \geq 5$, and suppose $E[p]$ is reducible. Then $E[p]^{\mathrm{ss}} \cong \mathbf{1} \oplus \bar\chi$, and some curve $E''$ isogenous to $E$ has a rational point of order $p$.
\end{proposition}

\begin{proof}
Let $0 \to \varphi \to E[p] \to \varphi' \to 0$ with characters $\varphi, \varphi' \colon G_\Q \to \F_p^\times$, so $\varphi\varphi' = \bar\chi$ by the Weil pairing. At a prime $\ell \mid N$ the curve has multiplicative reduction and inertia acts unipotently on $E[p]$, so both characters are unramified at $\ell$. At $p$ the reduction is good; it cannot be supersingular, since then $E[p]|_{G_p}$ would be irreducible (inertia acting through a fundamental character of level $2$), so it is ordinary, and on $G_p$ one of $\varphi, \varphi'$ is unramified and the other is $\bar\chi$ times an unramified character. Both characters are thus unramified outside $p$ and tamely ramified at $p$, so they are powers of $\bar\chi$, and the local structure at $p$ forces $\{\varphi, \varphi'\} = \{\mathbf{1}, \bar\chi\}$.

If $\varphi = \mathbf{1}$, take $E'' = E$. If $\varphi = \bar\chi$, the subgroup is isomorphic to $\mu_p$; put $E'' = E/\mu_p$. The image of $E[p]$ in $E''$ is a subgroup of order $p$ on which $G_\Q$ acts through $\varphi' = \mathbf{1}$, giving a rational point of order $p$.
\end{proof}

\begin{theorem}\label{thm:irred}
Let $p \geq 5$, and let $E/\Q$ have good reduction at $p$ and squarefree conductor $N$.
\begin{enumerate}
\item If $p \geq 11$, then $E[p]$ is absolutely irreducible.
\item If $p \in \{5, 7\}$ and $\ell \not\equiv \pm 1 \pmod p$ for every $\ell \mid N$, then $E[p]$ is absolutely irreducible.
\end{enumerate}
In particular, conditions (1), (2), (4) of \cite[Def.5.1]{RW} imply condition (3).
\end{theorem}

\begin{proof}
Since $\bar\rho$ is odd and $p$ is odd, complex conjugation has the two distinct eigenvalues $\pm 1$ in $\F_p$, so any $\bar\rho$-stable line over $\bar\F_p$ is one of its eigenlines, which are defined over $\F_p$; hence irreducibility implies absolute irreducibility. Suppose $E[p]$ is reducible, and let $E''$ be as in Proposition \ref{prop:shape}.

(1) By Mazur's theorem \cite{M}, a rational point of prime order on an elliptic curve over $\Q$ has order $2, 3, 5$ or $7$, so $p \leq 7$.

(2) Since isogenous curves have the same primes of bad reduction, every prime of bad reduction of $E''$ divides $N$, hence is $\neq p$ and $\not\equiv \pm 1 \pmod p$. By a theorem of Yasuda \cite[Thm.1.2]{Y12} (see \cite[Thm.3.3]{Y16} for the statement over $\Q$), for $p \in \{5,7\}$ an elliptic curve over $\Q$ whose primes of bad reduction are all $\not\equiv 0, \pm 1 \pmod p$ has no rational point of order $p$. This contradicts the existence of $E''$.
\end{proof}

\begin{remark}
Part (1) is classical. Yasuda's theorem does not require semistability; its proof in \cite{Y16} proceeds through the Tate curve and the vanishing of the class group of $\Q(\mu_p)$ for $p = 5, 7$. As a consistency check, among all $2{,}164{,}260$ isogeny classes of conductor at most $500{,}000$ and all $p \in \{5, 7, 11, 13, 17, 19, 23, 29\}$, no class satisfying (1), (2), (4) has reducible mod-$p$ image in Cremona's data (Section \ref{sec:repro}).
\end{remark}

\section{An upper bound for unobstructed curves}\label{sec:upper}

\subsection{Setup}
Every elliptic curve over $\Q$ is isomorphic to a unique $E_{a,b}\colon y^2 = x^3 + ax + b$ with $(a,b) \in \Z^2$, $4a^3 + 27b^2 \neq 0$, and $(a,b)$ \emph{minimal}, i.e.\ there is no prime $q$ with $q^4 \mid a$ and $q^6 \mid b$. We set $H(E_{a,b}) = \max(4|a|^3, 27 b^2)$ (the normalization $\max(|a|^3, b^2)$ of \cite{RW} changes nothing below) and
$$\B(X) = \{(a,b) \in \Z^2 : |a| \leq A,\ |b| \leq B\}, \qquad A = (X/4)^{1/3},\ B = (X/27)^{1/2},$$
so $|\B(X)| \asymp X^{5/6}$, and the number of elliptic curves of height at most $X$ is $\asymp X^{5/6}$. Write $\Delta(a,b) = -16(4a^3 + 27b^2)$.

For a prime $\ell \geq 5$ let
$$T_\ell = \{(a,b) \in (\Z/\ell)^2 :\ 4a^3 + 27b^2 \equiv 0,\ a \not\equiv 0\}.$$

\begin{lemma}\label{lem:Tl}
Let $\ell \geq 5$ be prime. Then $\#T_\ell = \ell - 1$, and for minimal $(a,b)$ the curve $E_{a,b}$ has multiplicative reduction at $\ell$ if and only if $(a,b) \bmod \ell \in T_\ell$.
\end{lemma}

\begin{proof}
For fixed $a \not\equiv 0$ the congruence $27 b^2 \equiv -4a^3$ has $1 + \left(\frac{-3a}{\ell}\right)$ solutions $b$, since $-4a^3/27$ and $-3a$ differ by a nonzero square. Summing over the $\ell - 1$ nonzero $a$, the Legendre symbols cancel, giving $\#T_\ell = \ell - 1$. For $\ell \geq 5$ we have $c_4 = -48a$ and $c_6 = -864b$, so a minimal pair gives a minimal model at $\ell$, and $E_{a,b}$ has multiplicative reduction at $\ell$ exactly when $\ell \mid \Delta$ and $\ell \nmid c_4$, i.e.\ when $(a,b) \bmod \ell \in T_\ell$.
\end{proof}

Note that $R$ contains neither $2$ nor $3$ when $p \geq 5$, so no local analysis at $2$ and $3$ is needed.

\subsection{The sieve}
For squarefree $d$ composed of primes $\ell \geq 5$, put $g(d) = \prod_{\ell \mid d} (\ell - 1)/\ell^2$ and
$$\mathcal{A}_d = \#\{(a,b) \in \B(X) :\ (a,b) \bmod \ell \in T_\ell \text{ for every } \ell \mid d\}.$$

\begin{lemma}\label{lem:level}
For such $d$, $\mathcal{A}_d = g(d)|\B(X)| + r_d$ with $|r_d| \ll A + B + d$, the implied constant being absolute.
\end{lemma}

\begin{proof}
By the Chinese Remainder Theorem and Lemma \ref{lem:Tl}, the condition defines $\prod_{\ell \mid d}(\ell - 1) \leq d$ residue classes of pairs modulo $d$. Each class meets $\B(X)$ in $(2A/d + \theta_1)(2B/d + \theta_2)$ points with $|\theta_i| \leq 1$, and $|\B(X)| = 4AB + O(A + B + 1)$. Summing over the classes gives the claim.
\end{proof}

\begin{theorem}\label{thm:upper}
Fix $p \geq 5$. Let $U_p(X)$ be the number of elliptic curves $E/\Q$ with $H(E) \leq X$, good reduction at $p$, $\bar\rho_{E,p}$ absolutely irreducible, and $\R_{E,p} \cong \Z_p$. Then
\begin{align*}
U_p(X) \leq M_p(X) &:= \#\{(a,b) \in \B(X) :\ (a,b) \bmod \ell \notin T_\ell \text{ for every } \ell \in R\} \\
&\ll_p X^{5/6}(\log X)^{-\kappa}.
\end{align*}
\end{theorem}

\begin{proof}
The first inequality follows from Proposition \ref{prop:obstruct} and Lemma \ref{lem:Tl}: a curve counted by $U_p(X)$ has no prime of multiplicative reduction in $R$, and its minimal pair lies in $\B(X)$.

For the second, let $z = X^{1/40}$, $P(z) = \prod_{\ell \in R, \ell < z} \ell$, and $D = z^{10} = X^{1/4}$. Then $M_p(X)$ is at most the number of pairs in $\B(X)$ avoiding $T_\ell$ for all $\ell \mid P(z)$. We use Mertens' theorem for a fixed modulus in the sharp form of Williams \cite{Wil}: for $(a,p) = 1$,
$$\prod_{\substack{\ell \leq x \\ \ell \equiv a \ (p)}}\Big(1 - \frac1\ell\Big) = \frac{C(p,a)}{(\log x)^{1/(p-1)}}\Big(1 + O_p\Big(\frac{1}{\log x}\Big)\Big), \qquad C(p,a) > 0.$$
Taking the product over $a = \pm 1$ gives $\prod_{\ell \in R, \ell \leq x}(1 - 1/\ell) = C_p (\log x)^{-\kappa}(1 + O_p(1/\log x))$ with $C_p > 0$. Since $(1 - g(\ell)) = (1 - 1/\ell)\big(1 + \frac{1}{\ell(\ell-1)}\big) \geq 1 - 1/\ell$, it follows that the density $g$ satisfies the dimension condition
$$\prod_{\substack{w \leq \ell < z \\ \ell \in R}}\big(1 - g(\ell)\big)^{-1} \leq \Big(\frac{\log z}{\log w}\Big)^{\kappa}\Big(1 + \frac{K_p}{\log w}\Big) \qquad (2 \leq w < z).$$
The fundamental lemma of sieve theory in dimension $\kappa$ (see \cite{HR} or \cite{FI}), with level $D$ and the remainders of Lemma \ref{lem:level}, gives
$$M_p(X) \ll_p |\B(X)| V(z) + \sum_{d \mid P(z),\ d < D} |r_d|, \qquad V(z) = \prod_{\ell \in R,\ \ell < z}\big(1 - g(\ell)\big).$$
By the same estimate, $V(z) = \prod_{\ell \in R, \ell < z}(1 - 1/\ell)\big(1 + \frac{1}{\ell(\ell-1)}\big) \ll_p (\log z)^{-\kappa} \ll_p (\log X)^{-\kappa}$, and the remainder is $\ll D(A + B + D) \ll X^{3/4}$.
\end{proof}

\begin{corollary}\label{cor:Ep}
The number of curves of height at most $X$ lying in $\Ep$ is $\ll_p X^{5/6}(\log X)^{-\kappa}$.
\end{corollary}

\begin{proof}
A curve in $\Ep$ has squarefree conductor, so all its bad primes are multiplicative, and none lies in $R$; its minimal pair is therefore counted by $M_p(X)$.
\end{proof}

\begin{remark}[Lower bounds]\label{rem:lower}
We expect $M_p(X) \asymp_p X^{5/6}(\log X)^{-\kappa}$, and likewise for the height-ordered count of $\Ep$. The integer analogue is classical: by the Selberg-Delange method (see \cite{T}), the number of squarefree $n \leq x$ with no prime factor in $R$ is $\asymp_p x(\log x)^{-\kappa}$. As a numerical check of the exponent, the proportion of squarefree integers coprime to $p$ with no prime factor in $R$ drops between $x = 4000$ and $x = 500{,}000$ by the factors $0.788$ and $0.859$ for $p = 5, 7$, against $(\log 4000/\log 500{,}000)^{\kappa} = 0.795, 0.858$. For the discriminant, the obstacle to a lower bound is the contribution of $(a,b)$ for which $\Delta(a,b)$ has a prime factor in $R$ larger than the level of distribution available from Lemma \ref{lem:level}; heuristically this contribution has the same order of magnitude as the main term, so controlling it requires sieve upper bounds for large prime factors of $\Delta(a,b)$ on sublattices. We do not pursue this here.
\end{remark}

\begin{conjecture}\label{conj:height}
For each $p \geq 5$, $U_p(X) \gg_p X^{5/6}(\log X)^{-\kappa}$.
\end{conjecture}

By Corollary \ref{cor:sandwich}, Conjecture \ref{conj:height} would follow from a lower bound of the same order for the height-ordered count of $\Ep$, together with the statement that a positive proportion of these curves have $p \nmid m_E$.

\section{Residual ramification and modular degrees}\label{sec:congruence}

Throughout this section $E$ runs over isogeny classes satisfying (1), (2), (4), which lie in $\Ep$ by Theorem \ref{thm:irred}. We write $\CL_p = 1 - \prod_{i \geq 1}(1 - p^{-i})$, so $\CL_5 = 0.239667\ldots$, $\CL_7 = 0.163205\ldots$, $\CL_{11} = 0.099167\ldots$, $\CL_{13} = 0.082838\ldots$.

\subsection{The residual-ramification locus}

\begin{proposition}\label{prop:lowering}
Let $E$ satisfy (1), (2), (4), and suppose $\bar\rho = \bar\rho_{E,p}$ is unramified at some $\ell \mid N$. Then $p \mid r_E$ and $p \mid m_E$.
\end{proposition}

\begin{proof}
Since $\ell \not\equiv 1 \pmod p$ and $\bar\rho$ is irreducible (Theorem \ref{thm:irred}), modular of level $N$ and unramified at $\ell \| N$, Mazur's principle (see \cite{R}) gives a newform $g$ of weight $2$ and level dividing $N/\ell$ with $\bar\rho_{g,\lambda} \cong \bar\rho$ for a prime $\lambda \mid p$ of its coefficient field. As $g$ has level a proper divisor of $N$, it is not Galois conjugate to $f_E$, so $p$ is a congruence prime for $f_E$ in the sense of \cite[Def.4.3]{RW}. By \cite[\S 5.1]{RW} the primes dividing $r_E$ are exactly these congruence primes, so $p \mid r_E$, and $p \mid m_E$ by \cite[Thm.2.1]{ARS} since $p \nmid N$.
\end{proof}

\begin{remark}
The implication from unramifiedness at $\ell$ to $p$ being a congruence prime is used in the proof of \cite[Thm.4.6]{RW}, where it is attributed to \cite[Lemma 12]{W2}. We record Proposition \ref{prop:lowering} only because it is the deterministic half of the decomposition below. In the data, the locus $L$ contains $21{,}809$, $25{,}544$, $16{,}973$ and $13{,}215$ classes of conductor at most $500{,}000$ for $p = 5, 7, 11, 13$, and $p \mid m_E$ in each of them.
\end{remark}

\begin{proposition}\label{prop:primevac}
Let $N = \ell$ be prime with $\ell \not\equiv 1 \pmod p$, and suppose $\bar\rho_{E,p}$ is irreducible. Then $\bar\rho_{E,p}$ is ramified at $\ell$, i.e.\ $p \nmid v_\ell(\Delta_E)$.
\end{proposition}

\begin{proof}
Otherwise Mazur's principle would produce $\bar\rho$ from $S_2(\mathrm{SL}_2(\Z)) = 0$.
\end{proof}

Consistently, among the $2{,}736$ (resp.\ $3{,}592$) prime-conductor classes of conductor at most $500{,}000$ satisfying (1), (2), (4) for $p = 5$ (resp.\ $7$), none lies in $L$.

\subsection{The decomposition}
For any finite set $\mathcal{F}$ of isogeny classes satisfying (1), (2), (4), write $P_{\mathcal{F}}$ for proportions within $\mathcal{F}$. Proposition \ref{prop:lowering} gives the exact identity
\begin{equation}\label{eq:decomp}
P_{\mathcal{F}}(p \mid m_E) = P_{\mathcal{F}}(L) + \big(1 - P_{\mathcal{F}}(L)\big) P_{\mathcal{F}}(p \mid m_E \mid E \notin L).
\end{equation}
Table \ref{tab:omega} shows both sides for $\mathcal{F}$ the classes of conductor at most $500{,}000$ stratified by $\omega(N)$. The unconditional frequency rises steadily with $\omega(N)$, nearly doubling for $p = 5$ between one and six bad primes, while the conditional frequency shows no statistically detectable variation across the observed strata. At $\omega(N) = 1$ the two agree, as Proposition \ref{prop:primevac} predicts.

\begin{table}[h]
\centering
\small
\begin{tabular}{llrrrrr}
\toprule
$p$ & $\omega(N)$ & classes & $P(p \mid m_E)$ & $\#(E \notin L)$ & $P(p \mid m_E \mid E \notin L)$ & s.e. \\
\midrule
$5$ & 1 & 2{,}736 & 24.01 & 2{,}736 & 24.01 & 0.82 \\
 & 2 & 17{,}815 & 32.60 & 15{,}830 & 24.15 & 0.34 \\
 & 3 & 37{,}712 & 39.71 & 29{,}870 & 23.88 & 0.25 \\
 & 4 & 31{,}354 & 45.41 & 22{,}568 & 24.16 & 0.28 \\
 & 5 & 8{,}973 & 49.76 & 5{,}889 & 23.45 & 0.55 \\
 & 6 & 293 & 51.19 & 181 & 20.99 & 3.03 \\
\midrule
$7$ & 1 & 3{,}592 & 16.23 & 3{,}592 & 16.23 & 0.62 \\
 & 2 & 26{,}555 & 21.47 & 24{,}890 & 16.22 & 0.23 \\
 & 3 & 65{,}490 & 25.66 & 58{,}096 & 16.20 & 0.15 \\
 & 4 & 68{,}479 & 29.49 & 57{,}614 & 16.20 & 0.15 \\
 & 5 & 27{,}437 & 32.30 & 22{,}178 & 16.25 & 0.25 \\
 & 6 & 1{,}946 & 32.12 & 1{,}585 & 16.66 & 0.94 \\
\midrule
$11$ & 1 & 4{,}316 & 9.24 & 4{,}316 & 9.24 & 0.44 \\
 & 2 & 34{,}635 & 12.08 & 33{,}761 & 9.81 & 0.16 \\
 & 3 & 97{,}692 & 13.72 & 93{,}386 & 9.74 & 0.10 \\
 & 4 & 120{,}804 & 14.97 & 113{,}742 & 9.69 & 0.09 \\
 & 5 & 60{,}160 & 15.67 & 56{,}008 & 9.42 & 0.12 \\
 & 6 & 8{,}229 & 15.80 & 7{,}650 & 9.42 & 0.33 \\
\midrule
$13$ & 1 & 4{,}617 & 8.32 & 4{,}617 & 8.32 & 0.41 \\
 & 2 & 37{,}518 & 9.38 & 36{,}946 & 7.98 & 0.14 \\
 & 3 & 110{,}362 & 10.82 & 107{,}262 & 8.24 & 0.08 \\
 & 4 & 145{,}157 & 11.59 & 139{,}672 & 8.12 & 0.07 \\
 & 5 & 81{,}031 & 12.14 & 77{,}532 & 8.17 & 0.10 \\
 & 6 & 12{,}113 & 12.23 & 11{,}554 & 7.98 & 0.25 \\
\bottomrule
\end{tabular}
\medskip
\caption{Frequencies of $p \mid m_E$ (in percent) among isogeny classes of conductor at most $500{,}000$ satisfying (1), (2), (4), by number of bad primes; s.e.\ is the binomial standard error $\sqrt{\hat q(1-\hat q)/n}$ of the conditional frequency. Numerators are recoverable from the table.}\label{tab:omega}
\end{table}

\subsection{The conditional frequency}

\begin{conjecture}\label{conj:CL}
Fix $p \geq 5$. Ordering isogeny classes satisfying (1), (2), (4) by conductor,
$$\lim_{X \to \infty} P_{N \leq X}\big(p \mid m_E \mid E \notin L\big) = \CL_p,$$
and for each $k \geq 1$ the same limit holds among classes with $\omega(N) = k$.
\end{conjecture}

The evidence at $X = 500{,}000$ is as follows.

\emph{Variation across $\omega(N)$.} A chi-square test of homogeneity of the conditional frequency across $\omega(N) \in \{1, \ldots, 6\}$ gives $p$-values $0.76$, $0.998$, $0.25$, $0.64$ for $p = 5, 7, 11, 13$. Thus there is no statistically detectable variation across the observed strata. This is weaker than independence of $\omega(N)$: a failure to reject homogeneity indicates only that any variation is not large compared with the standard errors in Table \ref{tab:omega}, which range from under a tenth of a percentage point in the largest strata to about three points in the smallest.

\emph{The value.} The pooled conditional frequencies are $23.98 \pm 0.15$, $16.21 \pm 0.09$, $9.65 \pm 0.05$ and $8.15 \pm 0.04$ percent, against $100\CL_p = 23.97$, $16.32$, $9.92$, $8.28$. For $p = 5, 7$ these agree within sampling error; for $p = 11, 13$ the pooled frequency is below $\CL_p$ by about five and three standard errors respectively. Table \ref{tab:windows} shows that the shortfall is concentrated at small conductors and shrinks across disjoint conductor windows; in the largest window it is $0.16 \pm 0.06$ for $p = 11$ and $0.02 \pm 0.05$ for $p = 13$. Table \ref{tab:cutoff} shows the cumulative frequencies rising toward $\CL_p$ as the cutoff grows. For $p = 5$ the cumulative frequency at $X = 500{,}000$ slightly exceeds $\CL_5$, and two windows give negative shortfalls, so the approach is not monotone from below; the data do not support any particular rate of convergence.

\begin{table}[h]
\centering
\small
\begin{tabular}{lcccc}
\toprule
window & $p=5$ & $p=7$ & $p=11$ & $p=13$ \\
\midrule
$N \leq 4000$ & $4.68 \pm 1.29$ & $4.45 \pm 0.82$ & $4.22 \pm 0.53$ & $3.94 \pm 0.45$ \\
$4000 < N \leq 20{,}000$ & $0.65 \pm 0.70$ & $1.20 \pm 0.43$ & $1.04 \pm 0.27$ & $0.26 \pm 0.22$ \\
$20{,}000 < N \leq 10^5$ & $-0.39 \pm 0.35$ & $0.41 \pm 0.21$ & $0.30 \pm 0.13$ & $0.38 \pm 0.11$ \\
$10^5 < N \leq 5 \cdot 10^5$ & $-0.05 \pm 0.18$ & $-0.10 \pm 0.10$ & $0.16 \pm 0.06$ & $0.02 \pm 0.05$ \\
\midrule
classes ($E \notin L$) & 1{,}094; 3{,}736; & 2{,}047; 7{,}524; & 3{,}211; 12{,}549; & 3{,}731; 15{,}036; \\
 & 14{,}921; 57{,}323 & 31{,}310; 127{,}074 & 55{,}816; 237{,}287 & 67{,}918; 290{,}898 \\
\bottomrule
\end{tabular}
\medskip
\caption{Shortfall $100\CL_p - P(p \mid m_E \mid E \notin L)$, in percentage points, in disjoint conductor windows, with the standard error under the hypothesis that the frequency equals $\CL_p$. The last rows give the number of classes not in $L$ in each window.}\label{tab:windows}
\end{table}

\begin{table}[h]
\centering
\small
\begin{tabular}{lcccc|c}
\toprule
 & $X = 4000$ & $20{,}000$ & $100{,}000$ & $500{,}000$ & $100\CL_p$ \\
\midrule
$p=5$ & $19.29 \pm 1.19$ & $22.40 \pm 0.60$ & $23.88 \pm 0.30$ & $23.98 \pm 0.15$ & 23.97 \\
$p=7$ & $11.87 \pm 0.71$ & $14.43 \pm 0.36$ & $15.56 \pm 0.18$ & $16.21 \pm 0.09$ & 16.32 \\
$p=11$ & $5.70 \pm 0.41$ & $8.23 \pm 0.22$ & $9.31 \pm 0.11$ & $9.65 \pm 0.05$ & 9.92 \\
$p=13$ & $4.34 \pm 0.33$ & $7.29 \pm 0.19$ & $7.77 \pm 0.09$ & $8.15 \pm 0.04$ & 8.28 \\
\midrule
$P(L)$, $p=5$ & 20.67 & 20.18 & 21.20 & 22.06 & \\
$P(L)$, $p=7$ & 10.84 & 12.08 & 12.70 & 13.20 & \\
$P(L)$, $p=11$ & 4.12 & 4.62 & 4.91 & 5.21 & \\
$P(L)$, $p=13$ & 2.43 & 2.78 & 3.14 & 3.38 & \\
\bottomrule
\end{tabular}
\medskip
\caption{Cumulative conditional frequency $P(p \mid m_E \mid E \notin L)$ for $N \leq X$, in percent with binomial standard errors, and the proportion $P(L)$ of the residual-ramification locus.}\label{tab:cutoff}
\end{table}

\subsection{A heuristic for Conjecture \ref{conj:CL}}\label{sec:heuristic}

We explain why one should expect the locus $L$ to carry the variation in $\omega(N)$, and the complementary frequency to equal the Cohen-Lenstra constant. The explanation is the Euler characteristic formula \eqref{eq:wiles}, read as a statement about the shape of a presentation.

For $E$ satisfying (1), (2), (4), write
$$\delta(E) = \sum_{\ell \mid N} \dim H^0\big(G_\ell, \mathrm{Ad}^0\bar\rho(1)\big)$$
for the right side of \eqref{eq:wiles}. Every $\ell \mid N$ satisfies $\ell \not\equiv \pm 1 \pmod p$, so by Lemma \ref{lem:local} the $\ell$-th term vanishes when $\bar\rho$ is ramified at $\ell$; when $\bar\rho$ is unramified at $\ell$ it is split, $\mathrm{Ad}^0\bar\rho(1)|_{G_\ell} \cong \bar\chi^2 \oplus \bar\chi \oplus \mathbf{1}$, and the term equals $1$. Hence $\delta(E)$ is the number of bad primes at which $\bar\rho$ is unramified, and
$$\delta(E) = 0 \iff E \notin L.$$

Equation \eqref{eq:wiles} says that the tangent space of $\R_{E,p}$ and its dual Selmer group differ in dimension by exactly $\delta(E)$. These are the kernel and cokernel dimensions of an $\F_p$-linear map $\F_p^{n+\delta(E)} \to \F_p^{n}$ for suitable $n$, and the Cohen-Lenstra philosophy \cite{CL} suggests modelling this map as a matrix chosen uniformly among all $p^{n(n+\delta(E))}$ of them, with $n$ large. This model predicts:
\begin{itemize}
\item if $\delta(E) \geq 1$, the kernel is always nonzero, in agreement with Propositions \ref{prop:obstruct} and \ref{prop:lowering};
\item if $\delta(E) = 0$, the matrix is square and the kernel is nonzero for a proportion
$$1 - \frac{|\mathrm{GL}_n(\F_p)|}{p^{n^2}} = 1 - \prod_{i=1}^{n}\big(1 - p^{-i}\big) \longrightarrow \CL_p \qquad (n \to \infty)$$
of the matrices.
\end{itemize}
The model sees $E$ only through $\delta(E)$. The number of bad primes enters \eqref{eq:wiles} only through local terms, and off $L$ all of these vanish, so the predicted frequency is the same in every stratum of fixed $\omega(N)$.

Two facts connect the tangent space to the modular degree off $L$. First, if $\ell \mid N$ and $\bar\rho$ is ramified at $\ell$, then $H^0(G_\ell, \mathrm{Ad}^0\bar\rho(1)) = 0$, so the local Euler characteristic formula and local duality give $\dim H^1(G_\ell, \mathrm{Ad}^0\bar\rho) = \dim H^0(G_\ell, \mathrm{Ad}^0\bar\rho) = \dim H^1_{\mathrm{unr}}(G_\ell, \mathrm{Ad}^0\bar\rho)$. Thus the absence of a local condition at $\ell$ coincides with the unramified condition, and for $E \notin L$ the tangent space is the mod-$p$ counterpart of the Bloch-Kato Selmer group $H^1_\emptyset(\Q, W_p(E))$ of \cite[\S 4]{RW}; the two can differ only through the condition at $p$ (Remark \ref{rem:converse}). Second, for $E \in \Ep \setminus L$, $p \mid m_E$ if and only if $p$ is a strict congruence prime for $f_E$ (Remark \ref{rem:converse}), which by \cite[Prop.4.4]{RW} holds if and only if $H^1_\emptyset(\Q, W_p(E)) \neq 0$. So off $L$ the event $p \mid m_E$ is the nonvanishing of a Selmer group whose Euler characteristic is $0$, and the model predicts frequency $\CL_p$; on $L$ the Euler characteristic is positive and the event is forced. This is the decomposition \eqref{eq:decomp}, with Watkins' constant attached to the balanced part only. At prime conductor $L$ is empty (Proposition \ref{prop:primevac}), so there the unconditional heuristic of \cite{Wat, RW} and the model agree.

Two remarks on the scope of the model. The Selmer group and its dual are attached to the different modules $\mathrm{Ad}^0\bar\rho$ and $\mathrm{Ad}^0\bar\rho(1)$, which are not isomorphic for $p \geq 5$ (their determinants are $\mathbf{1}$ and $\bar\chi^3$). So there is no pairing on a single group of the kind that alters the predicted distribution for Selmer groups of elliptic curves themselves, and the plain square-matrix count is the natural guess. On the other hand, the model does not account for the shortfall at small conductors in Table \ref{tab:windows}, which is of similar size for all four primes; truncating the model at a finite $n$ would instead give a shortfall $\prod_{i=1}^{n}(1 - p^{-i}) - \prod_{i \geq 1}(1 - p^{-i}) \asymp p^{-n-1}$, which depends strongly on $p$.

\begin{remark}
Congruences with rational newforms of the same level play almost no role. For $p = 5$ (resp.\ $7$), only $0.27\%$ (resp.\ $0.03\%$) of the classes have another isogeny class of the same conductor whose traces $a_q$, $q < 100$, $q \neq p$, agree modulo $p$ and whose Atkin-Lehner signs agree at the bad primes; in each of these $268$ (resp.\ $64$) cases $p \mid m_E$. This agreement of finitely many traces is a numerical filter rather than a proof of a congruence (a Sturm-bound certificate would need far more coefficients), and no result of this paper depends on it.
\end{remark}

\section{The size of the residual-ramification locus}\label{sec:Lstats}

Combining Corollary \ref{cor:sandwich}, Proposition \ref{prop:lowering} and Conjecture \ref{conj:CL}, one is led to expect that, in the conductor ordering, the proportion of $\Ep$ with $p \nmid m_E$ is asymptotic to $(1 - \CL_p)\big(1 - P_X(L)\big)$. Whether $1 - P_X(L)$ stays bounded away from $0$ is a genuine question, and the data do not settle it.

For a prime $\ell$ write $q_\ell(p)$ for the proportion of classes satisfying (1), (2), (4) with $\ell \mid N$ for which $p \mid v_\ell(\Delta_E)$. Table \ref{tab:qsmall} records these proportions for small $\ell$. The sum $\sum_\ell q_\ell(p)P(\ell \mid N)$ is the expected number of bad primes at which $\bar\rho$ is unramified. A summable bound $q_\ell(p)P(\ell \mid N) \leq b_\ell$ holding uniformly in the cutoff would control the contribution of large primes to $L$, since the proportion of classes in $L$ by virtue of some $\ell > y$ is at most $\sum_{\ell > y} b_\ell$. It would not by itself show that $\limsup_X P_X(L) < 1$: that requires, in addition, information about the joint distribution of the events at the small primes, since even two dependent events can together cover the whole family. Two features of the data caution against reading such a decay off a single cutoff.

\begin{table}[h]
\centering
\small
\begin{tabular}{lccccccccc}
\toprule
$\ell$ & 2 & 3 & 5 & 7 & 13 & 17 & 23 & 37 & 47 \\
\midrule
$q_\ell(5)$ & 15.11 & 12.64 & -- & 8.65 & 6.38 & 5.92 & 4.37 & 3.60 & 3.23 \\
$q_\ell(7)$ & 9.45 & 7.46 & 4.96 & -- & -- & 1.64 & 1.43 & 1.04 & 0.63 \\
\bottomrule
\end{tabular}
\medskip
\caption{$q_\ell(p)$ in percent at $X = 500{,}000$; `--' marks $\ell = p$ or $\ell \in R$. Every cell is based on at least $4{,}500$ classes; the counts are in Appendix \ref{app:counts}.}\label{tab:qsmall}
\end{table}

\emph{Dependence on the cutoff.} At a fixed range of $\ell$, $q_\ell(p)$ keeps increasing with $X$ (Table \ref{tab:qcut}). The apparent decay of $q_\ell(p)$ in $\ell$ at a fixed cutoff is therefore in part a finite-size effect.

\emph{Dependence on the cofactor.} If $\ell \| N$ and $\bar\rho$ is unramified at $\ell$, then by Mazur's principle $\bar\rho$ arises from a newform of level dividing $M = N/\ell$. Since $N \leq X$, large $\ell$ forces small $M$, and for squarefree $M < 11$ there are no weight-$2$ newforms of any level dividing $M$ (each such level $d$ has $X_0(d)$ of genus $0$). In the data, $q_\ell(p) = 0$ for every class with $M < 11$, and for $\ell$ above $X/11$ these are the only classes present. Table \ref{tab:qcof} separates the two variables for $p = 5$.

\begin{table}[h]
\centering
\small
\begin{tabular}{lcccc}
\toprule
$\ell$ & $X = 4000$ & $20{,}000$ & $100{,}000$ & $500{,}000$ \\
\midrule
\multicolumn{5}{l}{$p = 5$} \\
$[50, 100)$ & 2.28 (5/219) & 2.62 (23/878) & 2.43 (90/3711) & 2.73 (395/14470) \\
$[100, 200)$ & 0 (0/238) & 0.69 (7/1019) & 1.36 (58/4266) & 1.75 (282/16112) \\
$[200, 400)$ & 0 (0/179) & 0.25 (2/796) & 0.89 (27/3025) & 1.07 (129/12068) \\
$[400, 800)$ & 0 (0/202) & 0.43 (3/691) & 0.46 (12/2624) & 0.66 (70/10584) \\
$[800, 1600)$ & 0 (0/120) & 0 (0/602) & 0.12 (3/2433) & 0.36 (34/9410) \\
\midrule
\multicolumn{5}{l}{$p = 7$} \\
$[50, 100)$ & 0.46 (2/437) & 0.43 (9/2085) & 0.67 (62/9268) & 0.67 (264/39349) \\
$[100, 200)$ & 0.24 (1/423) & 0.15 (3/2016) & 0.24 (22/9024) & 0.30 (110/36888) \\
$[200, 400)$ & 0 (0/301) & 0 (0/1352) & 0.12 (7/5788) & 0.15 (36/24645) \\
$[400, 800)$ & 0 (0/278) & 0 (0/1312) & 0.04 (2/5249) & 0.06 (14/22467) \\
\bottomrule
\end{tabular}
\medskip
\caption{$q_\ell(p)$ in percent, pooled over $\ell$ in dyadic ranges, for growing cutoff $X$; in parentheses, the number of pairs (class, $\ell$) with $p \mid v_\ell(\Delta_E)$ over the number of pairs.}\label{tab:qcut}
\end{table}

\begin{table}[h]
\centering
\small
\begin{tabular}{lccccc}
\toprule
$\ell$ & $M < 11$ & $[11, 50)$ & $[50, 200)$ & $[200, 1000)$ & $[1000, 5000)$ \\
\midrule
$[100, 200)$ & 0 (0/102) & 0 (0/400) & 1.47 (12/814) & 1.58 (72/4564) & 1.94 (198/10232) \\
$[400, 800)$ & 0 (0/256) & 0.58 (5/863) & 0.46 (8/1757) & 0.73 (55/7573) & 1.48 (2/135) \\
$[1600, 3200)$ & 0 (0/651) & 0 (0/2542) & 0.16 (7/4503) & 0.34 (4/1178) & \\
$[6400, 12800)$ & 0 (0/1878) & 0.10 (6/6004) & 0 (0/350) & & \\
\bottomrule
\end{tabular}
\medskip
\caption{$q_\ell(5)$ in percent at $X = 500{,}000$, by range of $\ell$ and of the cofactor $M = N/\ell$.}\label{tab:qcof}
\end{table}

\begin{question}\label{q:L}
Fix $p \geq 5$. In the conductor ordering of isogeny classes satisfying (1), (2), (4), is $\limsup_{X \to \infty} P_{N \leq X}(L) < 1$? Equivalently, does excluding residual ramification at every bad prime cost only a bounded factor?
\end{question}

The cumulative values of $P(L)$ in Table \ref{tab:cutoff} grow slowly with $X$ for each $p$. We do not know whether they converge. An affirmative answer to Question \ref{q:L}, together with Conjecture \ref{conj:CL}, would give a positive lower density for $\{p \nmid m_E\}$ inside $\Ep$ in the conductor ordering. In the height ordering, by contrast, the condition $p \mid v_\ell(\Delta_E)$ at a multiplicative prime requires $\ell^p \mid \Delta(a,b)$ with $\ell \nmid a$, which for each $a$ allows at most two residue classes of $b$ modulo $\ell^p$; it therefore has density $O(\ell^{-p})$, and the expected number of such primes is bounded by a convergent sum. This is one concrete way in which the two orderings differ, and it is why Conjecture \ref{conj:height} is stated for heights.

Finally, the proportion of isogeny classes in $\Ep$ itself, which by Corollary \ref{cor:Ep} is $O((\log X)^{-\kappa})$ in the height ordering, decreases visibly in the conductor ordering (Table \ref{tab:Ep}); we do not prove any statement for the conductor ordering.

\begin{table}[h]
\centering
\small
\begin{tabular}{lcccccccc}
\toprule
$p$ & 5 & 7 & 11 & 13 & 17 & 19 & 23 & 29 \\
\midrule
in $\Ep$, $N \leq 4000$ & 10.3 & 17.2 & 25.1 & 28.6 & 30.4 & 30.9 & 32.1 & 33.3 \\
in $\Ep$, $p \nmid m_E$, $N \leq 4000$ & 6.6 & 13.5 & 22.7 & 26.7 & 29.2 & 30.0 & 31.5 & 33.0 \\
in $\Ep$, $N \leq 500{,}000$ & 4.6 & 8.9 & 15.1 & 18.1 & 20.0 & 20.3 & 21.5 & 22.7 \\
in $\Ep$, $p \nmid m_E$, $N \leq 500{,}000$ & 2.7 & 6.5 & 12.9 & 16.0 & 18.5 & 19.0 & 20.5 & 21.9 \\
\bottomrule
\end{tabular}
\medskip
\caption{Percentage of all isogeny classes of conductor at most $X$ lying in $\Ep$, and lying in $\Ep$ with $p \nmid m_E$. The rows for $N \leq 4000$ agree with the table of \cite[\S 5.2]{RW}, except that \cite{RW} reports $22.3$ in the second row for $p = 11$.}\label{tab:Ep}
\end{table}

\section{Data and computations}\label{sec:repro}

All statistics are computed from Cremona's database \texttt{ecdata} \cite{E} at the fixed commit \texttt{25cec5ec}, which contains $2{,}164{,}260$ isogeny classes of conductor at most $500{,}000$. We use the files \texttt{degphi} (the modular degree of the first curve of each class, together with its prime divisors and $a$-invariants), \texttt{galrep} (Sutherland's labels for the non-maximal mod-$p$ images) and \texttt{aplist} (the coefficients $a_q$ for $q < 100$). The following conventions matter.

\begin{itemize}
\item \emph{Irreducibility.} $E[p]$ is reducible exactly when the image label has type \texttt{B} (Borel) or \texttt{Cs} (split Cartan, contained in a Borel); non-maximal labels of other types are irreducible. In the range considered, exactly one class carries a \texttt{Cs} label at a prime $p \leq 29$ without a \texttt{B} label (the class 11a at $p = 5$), and it does not satisfy (4).
\item \emph{Atkin-Lehner signs.} At a bad prime $q$, \texttt{aplist} records the eigenvalue $w_q$ of the Atkin-Lehner involution, so that $a_q = -w_q$. The only use of these signs is the comparison of two classes of the same level in the remark at the end of Section \ref{sec:congruence}, which is insensitive to this convention.
\item \emph{Valuations.} $v_\ell(\Delta_E)$ is computed from the minimal model of the first curve in the class. For classes satisfying (1), (2), (4), the condition $p \mid v_\ell(\Delta_E)$ is equivalent to unramifiedness of $\bar\rho$ at $\ell$ (Lemma \ref{lem:local}) and is the same for all curves in the class (Section \ref{sec:local}).
\item \emph{Modular degrees.} $p \mid m_E$ is read from the set of prime divisors of the modular degree recorded in \texttt{degphi}.
\end{itemize}

The computation is a single pass over the database in standard-library Python and runs in a few minutes. It reproduces both columns of the conductor-$4000$ table of \cite[\S 5.2]{RW} (with the one exception noted under Table \ref{tab:Ep}), and every number in this paper.

\subsection*{Data availability}
All data used in this paper are publicly available in Cremona's database \texttt{ecdata} \cite{E}; we used commit \texttt{25cec5ec}. The Python code that produces every table in this paper is available from the author upon request. The raw counts behind the tables are given in Appendix \ref{app:counts} and in the tables themselves.

\subsection*{Acknowledgements}
The Python code used for the computations was developed with the assistance of Claude Opus 5.5 (Anthropic); the author verified the code and its output and is responsible for all results. Claude Opus 5.5 is also used for language clarity and grammar checks. Every sentence, mathematical claim and final piece of wording in this manuscript is the author's own, and the author takes full responsibility for its content.

\appendix
\section{Raw counts}\label{app:counts}

Table \ref{tab:raw-family} gives, for each $p$ and cutoff $X$, the number of isogeny classes in $\Ep$ and the number of these with $p \mid m_E$; the total number of isogeny classes of conductor at most $X = 4000$, $20{,}000$, $100{,}000$, $500{,}000$ is $13{,}352$, $81{,}217$, $437{,}226$, $2{,}164{,}260$. Table \ref{tab:raw-q} gives the counts behind Table \ref{tab:qsmall}. The counts behind Tables \ref{tab:omega} to \ref{tab:qcof} are given in those tables.

\begin{table}[h]
\centering
\small
\begin{tabular}{lcccc}
\toprule
$p$ & $X = 4000$ & $20{,}000$ & $100{,}000$ & $500{,}000$ \\
\midrule
5 & 1379; 496 & 6051; 2303 & 25066; 10032 & 98883; 40292 \\
7 & 2296; 492 & 10886; 2696 & 46826; 12307 & 193499; 52776 \\
11 & 3349; 321 & 16523; 2060 & 75273; 10360 & 325836; 46791 \\
13 & 3824; 255 & 19304; 1905 & 89498; 9546 & 390798; 43979 \\
17 & 4057; 156 & 21056; 1258 & 98455; 7053 & 433781; 32906 \\
19 & 4129; 122 & 21156; 1090 & 99051; 5955 & 439291; 28052 \\
23 & 4287; 83 & 22317; 787 & 104766; 4679 & 465761; 22917 \\
29 & 4452; 51 & 23192; 549 & 110275; 3555 & 491802; 17710 \\
\bottomrule
\end{tabular}
\medskip
\caption{Number of isogeny classes in $\Ep$ with $N \leq X$; number of these with $p \mid m_E$.}\label{tab:raw-family}
\end{table}

\begin{table}[h]
\centering
\small
\begin{tabular}{lcc}
\toprule
$\ell$ & $p=5$ & $p=7$ \\
\midrule
2 & 10241/67766 & 12494/132270 \\
3 & 7039/55681 & 7990/107081 \\
5 & -- & 3698/74542 \\
7 & 2622/30308 & -- \\
13 & 1067/16715 & -- \\
17 & 786/13275 & 399/24363 \\
23 & 408/9341 & 242/16926 \\
37 & 207/5744 & 110/10619 \\
47 & 147/4545 & 52/8221 \\
\bottomrule
\end{tabular}
\medskip
\caption{Counts behind Table \ref{tab:qsmall}: classes with $p \mid v_\ell(\Delta_E)$ over classes with $\ell \mid N$, among classes of conductor at most $500{,}000$ satisfying (1), (2), (4).}\label{tab:raw-q}
\end{table}

\end{document}